\documentclass[final,a4paper,12pt]{amsart}
\usepackage{amsmath,amsthm,amsfonts,amssymb}

\usepackage{ifpdf} 
\ifpdf       
\usepackage[pdftex, linktocpage]{hyperref}
\else
\usepackage[hypertex, linktocpage]{hyperref}
\fi
\usepackage{verbatim}

\DeclareMathOperator{\NN}{\mathbb N}

\newcommand{\ca}[1]{{\mathcal #1}}
\newcommand{\rest}[1]{\mkern-3mu\upharpoonright\mkern-5mu \raisebox{-.5ex}{\ensuremath{#1}}}
\newcommand{\w}{\widetilde}

\theoremstyle{plain}
\newtheorem{theorem}{Theorem}
\newtheorem{lemma}[theorem]{Lemma}
\newtheorem{proposition}[theorem]{Proposition}
\newtheorem{corollary}[theorem]{Corollary}
\newtheorem{conjecture}[theorem]{Conjecture}
\newtheorem{claim}{Claim}

\newtheorem{question}[theorem]{Question}
\newtheorem{fact}[theorem]{Fact}
\theoremstyle{definition}
\newtheorem{remark}[theorem]{Remark}
\newtheorem{definition}[theorem]{Definition}
\newtheorem{example}[theorem]{Example}
\newtheorem{exercise}[theorem]{Exercise}

\numberwithin{theorem}{section}
\numberwithin{equation}{section}

\newcommand{\bt}{\begin{theorem}}
\newcommand{\et}{\end{theorem}}

\newcommand{\bl}{\begin{lemma}}
\newcommand{\el}{\end{lemma}}

\newcommand{\bd}{\begin{definition}}
\newcommand{\ed}{\end{definition}}

\newcommand{\beq}{\begin{equation}}
\newcommand{\eeq}{\end{equation}}

\newcommand{\bexa}{\begin{example}}
\newcommand{\eexa}{\end{example}}

\newcommand{\bexe}{\begin{exercise}}
\newcommand{\eexe}{\end{exercise}}

\newcommand{\bfact}{\begin{fact}}
\newcommand{\efact}{\end{fact}}

\newcommand{\bprop}{\begin{proposition}}
\newcommand{\eprop}{\end{proposition}}

\newcommand{\bp}{\begin{proof}}
\newcommand{\ep}{\end{proof}}

\newcommand{\bc}{\begin{corollary}}
\newcommand{\ec}{\end{corollary}}

\newcommand{\bq}{\begin{question}}
\newcommand{\eq}{\end{question}}

\newcommand{\bcong}{\begin{conjecture}}
\newcommand{\econg}{\end{conjecture}}

\newcommand{\br}{\begin{remark}}
\newcommand{\er}{\end{remark}}

\newcommand{\dirlim}{\varinjlim}

\begin{document}

\title[Cohomology of groups in o-minimal structures]{Cohomology of groups in o-minimal structures: acyclicity of the infinitesimal subgroup} 
\author{Alessandro Berarducci}

\address{Universit\`a di Pisa, Dipartimento di Matematica, Largo Bruno Pontecorvo 5, 56127 Pisa, Italy}
\email{\tt berardu@dm.unipi.it}
\thanks{Partially supported by the project:
Geometr\'{\i}a Real (GEOR) DGICYT MTM2005-02865 (2006-08)}
\date{Nov. 15, 2007} 

\vspace{-2em}
\begin{abstract} By recent work on some conjectures of Pillay,  each definably compact group in a saturated o-minimal structure is an expansion of a compact Lie group by a torsion free normal divisible subgroup, called its infinitesimal subgroup.  We show that the infinitesimal subgroup is cohomologically acyclic. This implies that the functorial correspondence between definably compact groups and Lie groups preserves the cohomology.
\end{abstract} 

\maketitle 
\section{Introduction}
This note is a continuation of \cite{B:07} and settles the main conjecture of that paper: the infinitesimal subgroup of a definably compact group in an o-minimal expansion of a field is cohomologically acyclic. The proof uses recent work in \cite{HP:07} on the ``compact domination conjecture''. The compact domination conjecture has so far been proved in the abelian case but we can handle the non-abelian case by a suitable reduction. As a corollary of the acyclicity, we obtain a canonical isomorphism between the (o-minimal) cohomology of a definably compact group and the cohomology of the real Lie group canonically associated to it. This answers a question in \cite{B:07}. In the abelian and semisimple case the isomorphism problem was also addressed in the recent preprint \cite{EJP:07} without passing through the acyclicity of the infinitesimal subgroup. 

We recall some definitions and facts. In this paper $G$ is always assumed to be a definably compact group in an o-minimal expansion $\mathbb M$ of a field. It is convenient to adopt the model theoretic convention of working in a ``universal domain'', so we assume that $\mathbb M$ is $\kappa$-saturated and $\kappa$-strongly homogeneous with $\kappa$ larger than the cardinality of most of the sets we will be interested in (except $\mathbb M$ itself !). In particular the language $L$ of $\mathbb M$ is {\bf small} (i.e. of  cardinality $<\kappa$) and $G$ is actually definable over a small elementary submodel $M \prec \mathbb M$. 
We are interested in a certain canonical subgroup $G^{00}$ of $G$ which, except in trivial cases, will not be definable but only type-definable. Recall that a set is called {\bf type-definable} if it is the intersection of a small family of definable sets, namely a family of definable sets indexed by a set of cardinality $<\kappa$. 
\br \label{small} For those that find themselves unconfortable with universal domains, let us point out that the notion of smallness may perhaps be better understood if we identity a definable set with a functor which to each model associates a set. A small family of definable sets could then be defined as a family whose index set does not depend on the model. Thus for instance if we work over an ordered field, $\bigcap_{n\in \NN} [-1/n,1/n]$ is a small intersection (i.e. type-definable), but $\bigcap_{t>0} [-t,t]$ is not, since its index set $\{t \mid t >0\}$ becomes larger if we enlarge the model (the first expression defines the infinitesimal elements of the model, the second one defines the singleton $0$). Whenever we say that a certain property of a given type-definable set holds, we implicity mean that it holds in all sufficiently saturated models (or equivalently in the universal domain). 
\er 
If $H$ is a normal type-definable subgroup of $G$ of {\bf bounded  index} (i.e. index $<\kappa$), the quotient $G/H$ can be endowed with a natural topology (the ``logic topology'', not the quotient one) making it into a compact group \cite{P:04}. It turns out that although $G$ and $H$ depend on the model, $G/H$ does not, in the sense that if $N\succ M$ are sufficiently saturated, $G(N)/H(N)$ is canonicallly isomorphic to $G(M)/H(M)$ as explained in \cite{P:04}. Note that this would not be true if the index of $H$ were not bounded: for instance $H$ could be the trivial subgroup. 

By \cite{BOPP:05} $G$ has the DCC (non-existence of infinite descending chains) on type-definable subgroups of $G$ of bounded index. The {\bf infinitesimal subgroup} $G^{00}$ of $G$ is defined as the smallest such subgroup (which exists by the DCC or by \cite{S:05}). 
The subgroup $G^{00}$ is necessarily normal and divisible, and in \cite{BOPP:05} it is shown that the compact group  $\Gamma = G/G^{00}$ (with the logic topology) is actually a Lie group (similarly for $G/N$ for any $N\lhd G$ type-definable of bounded index). So a posteriori the index of $G^{00}$ is $\leq 2^{\aleph_0}$ and it can be shown (see \cite{S:05} or \cite[Prop. 6.1]{HPP:07}) that $G^{00}$ can be type-defined over the same parameters needed to define $G$. By definition a subset of $\Gamma$ is closed in the {\bf logic topology} if and only if its preimage in $G$ under the natural map $\pi \colon G\to \Gamma$ is type-definable. We always consider on $G$ the topology of \cite{P:88}, namely the unique topology making $G$ at the same time a topological group and a ``definable $M$-manifold'' with a finite atlas. It turns out that $G^{00}$ is an open subset of $G$, and the morphism $\pi \colon G \to \Gamma$ is continuous with respect to the above topologies. (Warning: since $G^{00}$ is open in $G$, the quotient topology on $\Gamma = G/G^{00}$ coincides with the discrete topology, and not with the one we are considering.) In \cite{HPP:07} it is proved that in the abelian case $G^{00}$ is torsion free, so $\Gamma$ and $G$ have isomorphic torsion sub-groups. By \cite{EO:04} the torsion subgroup of a definably compact abelian group $G$ is isomorphic to the torsion subgroup of a torus of dimension $n$, where $n = \dim_M(G)$ is the o-minimal dimension of $G$. It then follows that $\dim_M(G) = \dim(\Gamma)$, where the latter is the dimension of $\Gamma$ as a Lie group.  This equality continues to hold even in the non-abelian case \cite{HPP:07}. Also the fact that $G^{00}$ is torsion free continues to hold in the non-abelian case \cite{B:07}. An important ingredient of the results in \cite{HPP:07} is the notion of generic set already studied in \cite{PP:07} making use of some work of A. Dolich \cite{Do:04}. A definable subset $X\subset G$ is {\bf generic} if finitely many left translates of $X$ cover $G$  (equivalently right-translates). The non-generic sets form an ideal \cite[Prop. 4.2]{HPP:07}, namely if the union of two definable sets is generic one of the two is generic. Since $G^{00}$ has bounded index it is easy to see that every definable set containing $G^{00}$ is generic. The converse fails, as a generic set may be disjoint from $G^{00}$. However using the results in \cite{HPP:07} it was shown in \cite{B:07} that $G^{00} = \bigcap_{X\mbox{ generic }}XX^{-1}$. We will need this caracterization in the sequel. 

\section{Compact domination}
Using a result in \cite{OP:07} 
in \cite{HP:07} it was proved that if $G$ is abelian, then $G$ is {\bf compactly dominated} in the sense of \cite{HPP:07}. By definition $G$ is compactly dominated (by $\pi \colon G \to \Gamma$) if given a definable set $X\subset G$, for all points $y\in \Gamma$ outside a set of Haar measure zero, $\pi^{-1}(y)$ is either contained in $X$ or in its complement. Said in other words $m(\pi(X) \cap \pi(X^c)) = 0$ where $m$ is the Haar measure on $\Gamma$. In \cite{B:04} it was suggested that one could try to define a probability measure $\mu$ on the boolean algebra of the definable subsets of a definably compact group $G$ by $\mu(X) = m(\pi(X))$.  
The problem of verifying the (finite) additivity of $\mu$ amounts to show that for two disjoint definable sets $A,B\subset G$, $m(\pi(A) \cap \pi(B)) = 0$, which is again equivalent to compact domination. Still another equivalent form of compact domination (see \cite{HPP:07} and \cite{HP:07}), is that the image under $\pi \colon G \to \Gamma$ of a definable set $X\subset G$ with empty interior has Haar measure zero. The equivalence with the original formulation follows from the the fact that $G^{00}$ (hence any of its translates) is open in $G$ and definably connected (as shown in \cite{BOPP:05}), namely it cannot meet a definable set and its complement without meeting its frontier. The following Proposition is the only place where compact domination is used in this paper. 
\bprop If $G$ is compactly dominated and $X\subset G$ is a definable generic set, then some left-translate (equivalently right-translate) of $X$ contains $G^{00}$. \eprop
\bp Write $G = g_1X \cup \ldots \cup g_k X$. By compact domination for $y \in \Gamma$ outside of a set of Haar measure zero and for each $i\in \{1, \ldots, k\}$, $\pi^{-1}(y)$ is either contained in $g_iX$ or in its complement. Since it must meet some $g_iX$, it must be contained there. 
\ep

\bt \label{abelian}
Suppose that $G$ is abelian, or more generally compactly dominated. Then there is a decreasing sequence $A_0\supset A_1 \supset A_2 \supset \ldots$ of definable subsets of $G$ such that $\bigcap_{n\in \NN} A_n = G^{00}$ and, for all $n$, $A_n$ is definably homeomorphic to a cell (in the o-minimal sense). \et
\bp
By \cite[Lemma 2.2]{B:07} there is a decreasing sequence $X_0\supset X_1 \supset X_2 \supset \ldots$ of definable subsets of $G$ such that $\bigcap_{n\in \NN} X_n = G^{00}$. Since $G^{00}$ is a group, $G^{00} = (\bigcap_n X_n)(\bigcap_n X_n)^{-1} $ and by the saturation assumption on on the model this equals $\bigcap_n X_n X_n^{-1}$. Define $A_n$ inductively as follows. 

Case $n=0$:  Write $X_0$ as a finite union of cells. Since $X_0 \supset G^{00}$, $X_0$ is generic. So at least one of its cells is generic. Let $C \subset X_0$ be such a cell. By compact domination there is $g\in G$ such that $Cg \supset G^{00}$. Let $A_0 = Cg$. 

Case $n+1$:  Suppose $A_n \supset G^{00} = \bigcap_i X_i X_i^{-1}$ has already been defined. By saturation there is $m$ such that $X_m X_m^{-1} \subset A_n$. Write $X_m$ as a finite union of cells. At least one of these cells is generic. Let $C \subset X_m$ be such a cell. By compact domination there is $g\in G$ such that $Cg \supset G^{00}$. Since $e \in G^{00}$, $g^{-1}\in C \subset X_m$. So $Cg \subset X_mX_m^{-1} \subset A_n$. Let $A_{n+1} = Cg$. 

Note that in case $n+1$ we can arrange the construction so that $m\geq n$, so in particular $m$ tends to $\infty$ with $n$. It follows that $\bigcap_i A_i \subset \bigcap_n X_n X_n^{-1} = G^{00}$, and therefore $\bigcap_i A_i = G^{00}$.  
\ep

\bfact \label{simple} The conclusion of Theorem \ref{abelian} was already known in the case when $G$ is a (non-abelian) definably simple definably compact group \cite[Thm. 8.5]{B:07}. \efact

\section{Types}
We have seen that we may consider a definable set $X$ as a functor that given a model $M$ (containing the parameters over which $X$ is defined) yields a set $X(M) \subset M^n$ (or simply think of $X$ as a formula and $X(M)$ as the set it defines). As usual we omit $M$ when it is implicit or irrelevant.  Given a definable set $X$ let $\w X (M)$ be the set of types over $M$ containing $X$ (which we can identify with ultrafilters of $M$-definable subsets of $X$). As in \cite{P:88b} We equip $\w X (M)$ with the {\bf spectral topology}: a basic open subset of $\w X$ is a set of the form $\w U$ with $U$ a definable relatively open subset of $X$. With this topology $\w X (M)$ is a quasi-compact normal spectral space, in general not Hausdorff (see \cite{CC:83} and \cite{P:88b}). 
On $\w X (M)$ one has also the {\bf constructible topology} which is compact Hausdorff and has as basic open sets the sets of the form $\w U$ with $U$ a definable subset of $X$ not necessarily open. We will however be exclusively interested in the spectral topology. 
Like $X$, also $\w X$ can be considered as a functor, namely the functor that given a model $M$ (containing the parameters over which $X$ is defined) yields the set $\w X (M)$ of all types over $M$ containing $X$. Given a definable map $f \colon X \to Y$ over $M$, we have an induced map $\w f \colon \w X \to \w Y$ (decorated with $M$ if needed) defined by $\w f (\alpha) = \{Z \mid f^{-1}(Z) \in \alpha\}$, where
$Z$ ranges over the $M$-definable sets. As noted in \cite{EJP:06}, given a definable function $f\colon X \to Y$ and definable subsets $A \subset X$ and $B \subset Y$, we have $\w f (\w A) = \w {f(A)}$ and $\w f ^{-1}(\w B) = \w {f^{-1}(B)}$. From this one readily deduces that if $f \colon X \to Y$ is continuous then so is $\w f \colon \w X \to \w Y$ with respect to the spectral topology. 
Moreover it is easy to see that the above equalities commute with the operation of restricting the types to a smaller model. We will also need a technical result (Lemma \ref{fibers} below) which, given a definable function $f\colon A \to B$, allows us to characterize $\w f ^{-1}(\beta) = \bigcap \{f^{-1}(Z) \mid Z \in \beta \}$, where $\beta$ is a type in $B$, in terms of $f^{-1}(b)$, where $b$ is a realization of $\beta$. 
The special case of Lemma \ref{fibers} when $f$ is the projection from $X\times [a,b]$ to $X$ was proved in \cite[Claim 4.5]{EJP:06} and was used there to prove the invariance of o-minimal sheaf cohomology under definable homotopies. The general case can be proved along similar lines. 
\bl \label{fibers} \footnote{We thank A. Fornasiero for providing a telephone proof of this lemma on demand.} Let $f\colon A \to B$ be a definable continuous function over the model $M$. Let $\beta \in \w B (M)$ be a type over $M$ containing $B$ and let $b\models \beta$ be a realization in some elementary extension of $M$. Let $M\langle b \rangle \succ M$ be the prime model over $M\cup \{b\}$. Let  $r \colon \w A (M\langle b \rangle) \to  \w A (M)$  be the map which sends a type over $M\langle b \rangle$ containing $A$ to its restriction to $M$. It is easy to see that $r$ is continuous but in general it is not an open map. However: \[r \colon \w{f^{-1}(b)} (M\langle b \rangle) \approx \w f ^{-1} (\beta) (M), \] namely $r$ sends the set of types over $M\langle b \rangle$ containing $f^{-1}(b)$ homeomorphically onto $\w f ^{-1}(\beta) \subset \w A (M)$. (This holds both in the spectral topology and in the constructible topology.) \el
\bp
The fact that $\w f ^{-1} (\beta)$ is the restriction of $\w {f ^{-1} (b)}$ to $M$ is clear: 
$\w f ^{-1} (\beta) = \w f ^{-1} (tp(b/M)) = \w {f^{-1}(b)} \rest M = r(\w {f^{-1}(b)})$.

Let us prove that $r \colon \w {f^{-1}(b)} \to \w f ^{-1} (\beta)$ is injective. So let $\alpha_1$ and $\alpha_2$ be two distinct types over $M\langle b \rangle$ containing the definable set $f^{-1}(b)$. Let $a_1 \models \alpha_1$ and $a_2 \models \alpha_2$ be two realizations. Then there is definable set over $M\langle b \rangle$ containing $a_1$ and not $a_2$. Since every element of $M\langle b \rangle$ is of the form $h(b)$ for some $M$-definable function $h$, there is a formula $\psi(x,y)$ over $M$ such that $\models \phi(a_1, b)$ and $\not \models \psi (a_2,b)$. But $b = f(a_1) = f(a_2)$. So $\models \phi(a_1, f(a_1))$ and $\not \models \phi(a_2, f(a_2))$.  This shows that $a_1,a_2$ have a different type over $M$ and concludes the proof of injectivity. 

The continuity is easy: a basic open set of $\w f ^{-1}(\beta)$ is of the form $\w U \cap \w f ^{-1}(\beta)$ for some $M$-definable open set $U$. Its preimage under $r$ is $\w U \cap \w {f ^{-1}(b)}$. 

It remains to prove that $r \colon \w {f^{-1}(b)} \to \w f ^{-1} (\beta)$ is an open map. We need:

\begin{claim} 
Let $n= \dim(\beta)$. We can assume that $B$ is a cell of dimension $n$ in $M^n$, $A \subset M^{n+k}$, and $f \colon A \to B$ is the projection onto the first $n$ coordinates. 
\end{claim}
In fact let $A' = \{(f(x),x) \mid x\in A\}$. Then $f \colon A \to B$ can be factored through the definable homeomorphism $A \to A'$ sending $x$ to $(f(x),x)$ followed by the projection $A' \to B$ sending $(u,v)$ to $u$. So we can assume that $f$ is a coordinate projection. 
To prove the rest of the claim consider a definable set $D$ of minimal dimension $n$ in $\beta$. We can assume that $D$ is contained in $B$ and is a cell. Replacing $B$ with $D$ we can assume that $B$ is a cell of dimension $n$ in some $M^m$ with $m\geq n$. We can further assume that $m = n$ since every cell of dimension $n$ is definably homeomorphic through a coordinate projection to an open cell in $M^n$. The claim is thus proved. \medskip 

To finish the proof of the Lemma consider a basic open subset of $\w {f^{-1}(b)}$, namely a subset of $\w A (M\langle b \rangle)$ of the form $\w U \cap \w {f ^{-1}(b)}$ where $U$ is an $M\langle b\rangle$-definable  open subset of $A$. We must show that the restriction of $\w U \cap \w {f ^{-1}(b)}$ to $M$ is open in $\w f ^{-1}(\beta)$. To this aim it suffices to find an $M$-definable open subset $L' \subset A$ such that $\w L' \cap \w {f ^{-1}(b) } = \w U \cap \w {f ^{-1}(b) }$ (as sets of types over $M\langle b \rangle$). Indeed the restriction would then be $\w {L'} \cap \w f^{-1}(\beta)$  (over $M$), which is an open subset of $\w f ^{-1}(\beta)$. To find the desired set $L'$ we reason as follows. Since $U$ is defined over $M\langle b \rangle$, we can write $U = U_b$ for some $M$-definable family $\{U_x \mid x\}$ of definable sets (not necessarily open).  Let $W$ be the set of all $x$ such that $U_x \cap f^{-1}(x)$ is relatively open in $f^{-1}(x)$. Then $W$ is an $M$-definable set containing $b$. So $\dim W = n$ and we can assume without loss of generality that $W$ is open. Let $L = \bigcup_{x \in W} ( U_x \cap f^{-1}(x) )$. Then $L$ is an $M$-definable set and $L \cap f^{-1}(b) = U_b \cap f^{-1}(b)$.  Recall that, thanks to the claim, $f$ is assumed to be the projection $(x,y) \mapsto x$. So $L$ is a set which projects onto the open set $W$ via $f$ and intersects each fiber $f^{-1}(x)$, with $x\in W$, into the relatively open subset $U_x \cap f^{-1}(x)$. This is not yet sufficient to ensure that $L$ is open in $A$. Consider however a cell decomposition of $L$. Its projection over $W$ gives a cell decomposition of $W$. Since $\dim(\beta) = n = \dim (W)$,  there is an open cell $W' \subset W$ of this decomposition containing $b$.  Now let $L' = L \cap f^{-1} (W')$. Then $L'$ is a cylinder over the open cell $W'$ consisting of a union of cells with common base $W'$. Since moreover $L'$ intersects each fiber $f^{-1}(x)$, with $x\in W'$, into the relatively open subset $U_x \cap f^{-1}(x)$, it follows now that $L'$ is open in $A$. We have thus found an $M$-definable open subset $L' \subset A$ such that $\w L' \cap \w {f ^{-1}(b) } = \w U \cap \w {f ^{-1}(b) }$ and this suffices to conclude.    
\ep 

Similarly to what already said for definable sets, also a type-definable set $Y$ can be considered as a kind of functor that to each model $M$ (containing the relevant parameters) associates a set $Y(M)$. Namely if $\{Y_i \mid i \in I\}$ is a small family of definable sets and $Y = \bigcap_{i\in I}Y_i$, then $Y(M) = \bigcap_{i\in I}Y_i(M)$. We can then define the space of types $\w Y (M)$ as $\bigcap_{i\in I} \w {Y_i}(M)$, and it easy to see that if $M$ is sufficiently saturated (with respect to the number of parameters in $Y$) then this is well defined, namely it does not depend on the representation of $Y(M)$ as an intersection. In other works we are saying that we can define
\[\w {\bigcap_{i\in I} Y_i} = \bigcap_{i\in I} \w {Y_i}\] (Meaning that the equality holds in sufficiently saturated models.) Note that the assumption that $I$ is small is essential, for instance considering the example in Remark \ref{small}, we have:
\[ \w {\{0\} } = \w {\bigcap_{t>0} [-t,t]} \subset \bigcap_{t >0} \w {[-t,t]} \] but the equality fails (in sufficienly saturated models). 

\bexa Although $G^{00}$ is open in $G$, it turns out that $\w {G^{00}}$ is closed in $\w G$ (and not open unless $G$ is trivial).  Indeed by \cite[Lemma 2.2]{B:07} $G^{00}$ can be written as a countable intersection $\bigcap_{n\in \NN} X_n$ of closed definable sets $X_n$, so $\w {G^{00}} = \bigcap_{n\in \NN} \w {X_n}$, which is closed. The fact that $\w {G^{00}}$ is not open in $\w G$ (except in trivial cases), follows from the fact that $G^{00}$ lives in the definably connected component $G^0$ of $G$ and if $X$ is definably connected then $\w X$ is connected. 
\eexa

\br \label{fiber2} Lemma \ref{fibers} continues to hold if $f \colon A \to B$ is the restriction of a definable function $g \colon X \to Y$ to two type-definable subsets $A\subset X$ and $B \subset Y$. This follows from the fact that for each $b\in B$, $f^{-1}(b) = g^{-1}(b) \cap A$, and $\w {g^{-1}(b) \cap A} = \w {g^{-1}(b)} \cap \w A$. \er

\section{Cohomology}
A sheaf cohomology theory for definable sets in o-minimal expansions $M$ of a group was studied in \cite{EJP:06} and \cite{BF:07}.  Given a definable set $X$ over $M$ and a sheaf $\ca F$ over $\w X (M)$, the cohomology $H^*(X; \ca F )$ is defined as the cohomology of the sheaf $\ca F$. Since $\w X (M)$ is a spectral space, it also coincides with the \v{C}ech cohomology of $\ca F$ \cite{CC:83}. When $Z$ is a fixed constant group of coefficients we write $H^*(\w X)$, or $H^*(X)$, for the cohomology $H^*(\w X (M))$ of the constant sheaf on $\w X (M)$ with stalk $Z$. A priori $H^*(X)$ depends on $M$, but if $M$ expands a field one can use the triangulation theorem and the results in \cite{EJP:06} to show that $H^*(X)$ does not actually depend on $M$, namely if $N\succ M$ then the restriction $r \colon \w X (N) \to \w X (M)$ induces is a canonical isomorphism $H^*(\w X(M)) \cong H^*(\w X(N))$. 

The cohomology of a type-definable set $Y$ can be defined similarly, namely $H^*(Y) = H^*(\w Y (M))$ where $\w Y (M)$ is defined as in the previous section (for $M$ sufficiently saturated). To compute the cohomology of a type-definable set and prove its invariance under elementary extensions one can use the following result. 
\bfact \label{limit} Let $X$ be a type-definable set, written as an intersection $\bigcap_i X_i$ of a small directed family of definable sets $X_i$. Then  ${H}^* (\widetilde {X}) = \dirlim_i {H}^*(\w {X_i})$. \efact
\bp This is a special case of \cite[Thm. 10.6]{Br:97}. See also \cite[Lemma C.3]{BF:07}. One needs to  observe that $\w X$ and $\w X_i$ are quasi-compact subsets of a spectral space and therefore they are ``taut'' subsets. This argument was used in \cite[Thm. 3.1]{De:85} and in \cite[Prop. 5.3.1]{J:06}, but  note that here we do not assume the $X_i$'s to be open. \ep 

\br If $M$ is only assumed to expand a group Fact \ref{limit} still holds. However for a definable set $X$ in an o-minimal expansion of a group the invariance of $H^*(\w X (M))$ under elementary extensions $N\succ M$ has so far been proved only when $X$ is definably compact \cite{BF:07}. Clearly it then also holds for those type-definable sets, such as $G^{00}$, which are small intersections of definably compact sets.\er 

\section{Acyclicity}
As usual let $G$ be a definably compact group in a sufficiently saturated o-minimal expansion $M$ of a field. 

\bt \label{acyclic} $\widetilde {G^{00}}$ is acyclic, namely $H^*(\w {G^{00}})$ is isomorphic to  the cohomology of a point (over any sufficiently saturated model).  
\et 
\bp By Fact \ref{simple} and Theorem \ref{abelian}, if $G$ is definably simple or abelian then $\w {G^{00}}$ is a decreasing intersection of a countable sequence of sets definably homeomorphic to cells. So by Fact \ref{limit} $\w {G^{00}}$ is acyclic in these case. (Similarly if $G$ is compactly dominated.) The general case can be reduced to the abelian and definably simple case as follows. Let 
\[ 0   \longrightarrow   H   \stackrel{i}{\longrightarrow}   G   \stackrel{\pi} 
{\longrightarrow}   B   \longrightarrow   0 \]
be an exact sequence of definable groups and definable group homomorphisms. Then by \cite[Thm 5.2]{B:07}
we obtain an induced sequence of type-definable groups:
\[ 0   \longrightarrow   H^{00}   \stackrel{i}\longrightarrow   G^{00} 
  \stackrel{\pi}\longrightarrow   B^{00}   \longrightarrow   0 \]
Passing to the space of types (over some model) we obtain a sequence of topological spaces and continuous maps 
\[ \w { H^{00}}   \stackrel{\w i }\longrightarrow   \w {G^{00}} 
  \stackrel{\w \pi }\longrightarrow  \w{ B^{00}}  \] 
It does not makes sense to ask whether the sequence is exact since these spaces do not carry a group structure. However, assuming that $\widetilde H^{00}$ is acyclic (i.e. $\w {H^{00}(M)}$ is acyclic for all $M$ sufficiently saturated), it follows that the fibers of $\w \pi$ are acyclic. Indeed if we work over the model $M$ and $\beta \in \w {B^{00}}(M)$,  then by Lemma \ref{fibers} the fiber $\w {\pi}^{-1}(\beta)$ is homeomorphic to $\w {H^{00}}(M\langle b \rangle)$ where $M\langle b \rangle$ is the prime model of a realization $b\models \beta$. (There is a small point to consider here. A priori we have defined $ {H^{00}}(N)$ only for $N$ sufficiently saturated, and $M\langle b \rangle$ need not be so even if $M$ is such. So what we mean is: let $\{X_i \mid i \in I\}$ be a small family of definable set such that $\w {H^{00}}(M) = \bigcap_i \w {X_i}(M)$ and define $\w {H^{00}} (M\langle b \rangle) = \bigcap_{i\in I} \w {X_i} (M\langle b \rangle)$. Now go to a saturated $N\succ M\langle b \rangle$, use the acyclicity assumption in $N$, and go back to $M\langle b \rangle$ using the invariance of cohomology under elementary extensions.)
Now by the appropriate version of the Vietoris-Begle theorem (see \cite[\S II, Thm. 11.7]{Br:97} and \cite[Thm. 4.3]{EJP:06}) under very general hypothesis a surjective continuous closed map with acyclic fibers induces an isomorphism in cohomology (with constant sheaves).  This is well known for maps between Hausdorff compact spaces, but it continues to hold for arbitrary topological spaces under the hypothesis that the fibers are ``taut''. This hypothesis is verified by $\w \pi$ since its fibers are type-definable, hence quasi-compact (see \cite[Prop. 2.20]{EJP:06}) in the spectral topology. Therefore we have an isomorphism $\w {\pi}^* \colon H^*(\widetilde {B^{00}}) \to H^*(\widetilde {G^{00}})$ over any sufficiently saturated model.  

It follows that the property of having an acyclic infinitesimal
subgroup is preserved under group extensions (i.e. if $0 \to H \to G \to B \to 0$ is exact and the property holds for $H$ and $B$ then it holds for $G$), 
and under isogenies (i.e. if it holds for $G$ and $H$ is finite, it holds for $B$). 
Since the property holds for the definably compact abelian groups, the definably compact definably simple groups, and obviously the finite groups, one gets it for the class $\ca P$ generated by these groups by extensions and isogenies, namely for every definably compact group. In fact suppose for a contradiction that $G$ is a definably compact groups of minimal dimension not belonging to $\ca P$. We can assume that $G$ is definably connected since the definably connected component $G^0$ of a definable group $G$ has finite index (so $G^0 \in \ca P$ iff $G \in \ca P$). If $G$ has an infinite definable abelian subgroup $N$, then $G/N$ has  lower dimension then $G$ and by the exact sequence $0 \to N \to G \to G/N \to 0$ we obtain $G$ as an extension of two groups in $\ca P$. So $G$ has no infinite definable abelian subgroups, namely $G$ is {\bf definably semisimple}. By \cite[Thm. 2.38]{PPS:00} if $G$ is definably connected and definably semisimple, then there are finitely many ``definably almost simple'' normal subgroups $H_1, \ldots, H_k$
of $G$ such that $G = \Pi_i H_i$ and  $H_i \cap \Pi_{j\neq i}H_j < Z(G)$, where $Z(G)$ is the finite center of $G$. Recall that a group $H$ is said to be {\bf definably almost simple} if it is non-abelian and has no infinite normal proper subgroup (we obtain an equivalent definition replacing ``infinite'' by ``definably connected and non-trivial'').  On the other hand if $H$ is definably almost simple, then its center $Z(H)$ is finite and $H/Z(H)$ is definably simple (using the fact that a finite normal subgroup $A$ of a definably connected group $H$ must be contained in its center). So all the $H_i$ belong to $\ca P$. Now let $R_m = \Pi_{i\leq m}H_i$. Considering the exact sequence $0 \to R_m \to R_m H_{m+1} \to R_m H_{m+1} /
R_m  \cong H_{m+1}/(R_m \cap H_{m+1} ) \to 0$ and noting that $R_m \cap H_{m+1}$ is finite (being contained in $Z(G)$), it follows by induction that each $R_m$, hence $G$, belongs to $\ca P$. 
\ep 

From the acyclicity of $G^{00}$ it follows, as explained in \cite[Cor. 8.8]{B:07}, that there is a natural isomomorphism in cohomology $H^*(G/G^{00}) \cong H^*(G)$, where the latter is defined as ${H}^*(\widetilde G )$. 
In fact let $\Psi \colon \widetilde G \to G/G^{00}$ be the map sending a type $\gamma \in \w G (M)$ ($M$ sufficiently saturated) to the coset $gG^{00}$ of a realization $g\models \gamma$ in some bigger model. Then $\Psi$ is a closed continuous map (independent of the choice of $g$) with acyclic fibers (as $\Psi^{-1}(\gamma) \approx \w {G^{00}} (M\langle g \rangle)$). So by the Vietoris-Begle theorem we obtain:

\bc 
$\Psi$ induces an isomorphism \[\Psi^* \colon H^*(G/G^{00}) \cong H^*(\w G)\] in cohomology. 
\ec

\subsection*{Acknowledgement} I thank Antongiulio Fornasiero for his comments on a preliminary version of this paper and for the proof of Lemma \ref{fibers}.

\end{document}